\documentstyle[12pt]{article} 
\newcommand{\doublespace}{
   \renewcommand{\baselinestretch}{1.2}
   \large\normalsize}

\def \Z{\Bbb Z}
\def \C{\Bbb C}

\def \Q{\Bbb Q}

\def \U{\cal U}
\def \S{\cal S}

\def \A{\cal A}
\def \M{\cal M}
\def \O{{\cal O}}
\def \wt{{\rm wt}}

\def \Res{{\rm Res}}

\def \End{{\rm End}}
\def \Hom{{\rm Hom}}

\def \Ind {{\rm Ind}}

\def \<{\langle} 
\def \>{\rangle} 

\def \a{\alpha }

\def \l{\lambda }
\def \L{\Lambda }

\def \b{\beta }
\def \om{\omega }

\def \ch{\chi}
\def \p{\phi}

\def \qed{\mbox{ $\square$}}
\def \pf {\noindent {\bf Proof:} \,}
\def \cg{\chi_g}
\def \cg'{\chi'_g}

\def \o{\otimes}
\def \d{\delta}

\def \r{\rho}

\input amssym.def
\input amssym
\doublespace
\begin{document}
\newtheorem{th}{Theorem}[section]
\newtheorem{thmn}{Theorem}
\newtheorem{prop}[th]{Proposition}
\newtheorem{cor}[th]{Corollary}
\newtheorem{lem}[th]{Lemma}
\newtheorem{rem}[th]{Remark}
\newtheorem{de}[th]{Definition}
\newtheorem{hy}[th]{Hypothesis}
\begin{center}
{\Large {\bf Vertex operator algebras, generalized doubles and dual pairs}
} \\
\vspace{0.5cm}
Chongying Dong\footnote{Supported by NSF grant 
DMS-9700923 and a research grant from the Committee on Research, UC Santa Cruz.} and Gaywalee Yamskulna\footnote{Supported by DPST grant.}
\\
Department of Mathematics, University of
California, Santa Cruz, CA 95064
\end{center}

\hspace{1.5 cm}
\begin{abstract} Let $V$ be a simple vertex operator algebra
and $G$ a finite automorphism group. 
Then there is a natural
right $G$-action on the set of all inequivalent irreducible $V$-modules.
Let $\S$ be a finite set of inequivalent irreducible $V$-modules which is 
closed under the action of $G.$ 
There is a finite dimensional
semisimple associative algebra $A_{\alpha}(G,\S)$ for
a suitable 2-cocycle naturally determined by the $G$-action on $\S$ 
such that $A_{\alpha}(G,\S)$ and the vertex operator algebra $V^G$ 
form a dual pair on the sum of $V$-modules in $\S$ in the sense of Howe.
In particular, every irreducible $V$-module is completely
reducible $V^G$-module.
\end{abstract}

\section{Introduction}

Let $V$ be a simple vertex operator algebra (see [B], [FLM]), and $G$ a finite 
automorphism group. A major problem in orbifold conformal 
field theory is to understand the module category
for the vertex operator algebra $V^G$ of $G$-invariants. 
In the case $V$ is holomorphic, this is related to the
quantum double (see [DPR] and [DM5]). 
The main feature in the study of orbifold theory is the appearance
of dual pairs. It is proved in [DLM1] and [DM2] that $G$ and
$V^G$ form a dual pair in the sense of Howe [H1]-[H2]. More precisely,
it is shown in [DLM1] that all the irreducible $G$-modules
occur in $V$ and the space of multiplicity of each irreducible
$G$-module in $V$ is an irreducible $V^G$-module. Moreover,
inequivalent irreducible $G$-modules produce inequivalent
$V^G$-modules in this way. 

In this paper we extend the duality result in [DLM1]
to any irreducible $V$-module and obtain again several duality
theorems of Schur-Weyl type. 
In the process we realize that it is better to consider 
a finite set of inequivalent irreducible modules which is $G$-stable
instead of one single module. The general setting is more natural
and the results are more beautiful. 

More explicitly, for an irreducible $V$-module $M=(M,Y_M)$ 
(see Section 4 for the details of the definition of a module) and $g\in G$
we define a new irreducible $V$-module $M\circ g=(M\circ g,Y_{M\circ g}).$ 
Here $M\circ g$ is equal to $M$ as a vector space and $Y_{M\circ g}(v,z)=
Y_M(gv,z)$ for $v\in V$ following [DM1]. A set $\S$ of irreducible
$V$-modules is called $G$-stable if for any $M\in S$ and $g\in G$ there
exists $N\in \S$ such that $M\circ g\cong N.$ Then every 
irreducible $V$-module $M$ produces naturally such a set 
by collecting all $M\circ g$ for $g\in G.$ Now we take 
a finite $G$-stable set $\S$ consisting of inequivalent 
irreducible $V$-modules. We construct a finite dimensional
semisimple associative algebra $A_{\alpha}(G,\S)$ for a
suitable 2-cocycle on $G$ (which is uniquely determined
up to isomorphism by $V,$ $G$ and $\S$) such that
$(A_{\alpha}(G,{\S}),V^G)$ forms a dual pair on ${\M}=\sum_{M\in \S}M$
in the precise sense of Howe [H1]-[H2]. That is, each simple
$A_{\alpha}(G,\S)$ occurs in $\M$ and its multiplicity space
is an irreducible $V^G$-module. Moreover, the different
multiplicity spaces are inequivalent $V^G$-modules.

These duality results not only tell us the complete
reducibility of  every irreducible $V$-module
as a $V^G$-module, but also provide 
an equivalence between the $A_{\alpha}(G,\S)$-module category
and a subcategory of $V^G$-modules generated by the irreducible
submodules of $\M$ by sending each simple $A_{\alpha}(G,\S)$-module
to its multiplicity space. This kind of idea has appeared in [DPR] and
[DM5] in the study of holomorphic orbifold conformal
field theory. 

Although the algebra $A_{\alpha}(G,\S)$ mentioned above appears
naturally in the theory of vertex operator algebras, the construction
itself is totally canonical and abstract. In fact, one can define
$A_{\alpha}(G,\S)$ for any finite group $G,$ any finite right $G$-set
$\S$ 
and a suitable 2-cocycle on $G$ (see Section 3). This algebra 
is essentially the crossed product in the theory 
of Hopf algebra (see [S], [BCM] and [DT]). It turns out that
$A_{\alpha}(G,\S)$ is the right algebra in the study of general
orbifold theory.  In the case that $\S$ is the dual basis of
$\C[G]^*,$ $A_{\alpha}(G,\S)$ is exactly the twisted double
$D_{\alpha}(G)$ introduced in [M] and [DM5] in the study of
holomorphic orbifold theory.  The twisted double $D_{\alpha}(G)$ is
conjecturally isomorphic to the twisted quantum double introduced in
[DPR] and [D].
  
A main tool in the proof of the main theorems is a series of
associative algebras $A_n(V)$ constructed in [DLM4] for 
nonnegative integers $n.$ The original motivation for introducing and
studying the $A_n(V)$ comes from the representation theory of vertex
operator algebras. Let $M=\sum_{n\geq 0}M(n)$ be an admissible
$V$-module with $M(0)\ne 0$ (see Section 4 for the definition
of admissible module). 
 The algebra $A_n(V)$ is a suitable
quotient of $V$ and ``takes care'' of the first $n+1$ pieces of $M$: each
$M(m)$ is a module for $A_n(V).$ In the case $n=0,$ $A_0(V)$ reduces
to the associative algebra $A(V)$ introduced previously in [Z]. The
main results concerning $A_n(V)$ are summarized in Section 4. The
result that we often use in this paper is the following: $M$ is
irreducible if and only if each $M(n)$ is a simple $A_n(V)$-module for
all $n\geq 0$ [DM4].  This key result allows us to reduce an infinite
dimensional problem ($M$ is infinite dimensional) to a finite
dimensional problem (each $M(n)$ is finite dimensional in our main
theorems).

This paper is organized as follows. In Section 2 we review twisted
group algebras. Section 3 is about the algebra $A_{\alpha}(G,\S)$
based on a finite group $G,$ a right $G$-set $\S$ and a suitable 
2-cocycle $\alpha$ on $G.$ We construct all simple modules for 
$A_{\alpha}(G,\S)$ explicitly by using twisted group
algebras. A basis for the center is also given for the future
study. In Section 4 we recall the various notions of modules for
a vertex operator algebra and review the theory on associative 
algebras $A_n(V).$ Section 5 is the first part on duality.
We show that if $M$ is an irreducible $V$-module which is $G$-stable
in the sense that $M\circ g\cong M$ for all $g,$ then
there exists a 2-cocycle $\a_M\in Z^2(G,\C^*)$ such that
$(\C^{\alpha_M}[G],V^G)$ forms a dual pair on $M$ where  $\C^{\alpha_M}[G]$
is the twisted group algebra. Section 6 is a continuation
of Section 5 without assuming that $M$ is $G$-stable. The algebra
$A_{\alpha}(G,S)$ based on a $G$-stable set $\S$ of $V$-modules 
enters the picture naturally. The two main theorems of this paper are proved in this
section. 

We would like to thank A. Wassermann for pointing out a gap
in the proof of Theorem 2.4 (i) of [DLM1]. This gap has been filled
in this paper (see Remark \ref{rcorect}).

\section{Twisted group algebras}
 
In this elementary section we review some results on twisted group
algebras which will be used in Section 3. 

Throughout this paper, $F^*$ denotes the multiplicative group of a field $F$ and $G$ a finite group which acts trivially on $F^*$. Let $\a$ be any element in $Z^2(G,F^*)$ . Thus $\a$ is a map 
$$\a:G\times G \rightarrow F^*$$
which satisfies the following properties:
\begin{eqnarray*}
& &\a(x,1)=\a(1,x)=1\ \ {\mbox for\ \ all}\ \ x\in G\\
& &\a(x,y)\a(xy,z)=\a(y,z)\a(x,yz)\ \ {\mbox for\ \ all}\ \ x,y,z\in G.
\end{eqnarray*}

\begin{de}
The {\em twisted group algebra} $F^{\a}[G]$ of $G$ over $F$ is defined
to be the vector space over $F$ with basis $\{\bar{g}|g\in
G\}$. Multiplication in $F^{\a}[G]$ is defined distributively by using
$$(a \bar{x})(b\bar{y})=ab\a(x,y)\overline{xy}\ \ (a,b\in F,x,y\in G).$$
\end{de}

It is easy to verify that $F^{\a}[G]$ is an $F$-algebra with $\bar{1}$ as the identity element and with 
$$\bar{g}^{-1}=\a(g^{-1},g)^{-1}\overline{g^{-1}}=\a(g,g^{-1})^{-1}\overline{g^{-1}}\ \ {\mbox for\ \ all}\ \ g\in G.$$  
\begin{lem}{[K]}
If {\sl char} $F$ does not divide $|G|$, then $F^{\a}[G]$ is semisimple.
\end{lem}

Our aim for this section is to find an explicit basis for the center $Z(F^{\a}[G])$.
\begin{de}
An element $g\in G$ is said to be $\a-{\em regular}$ if 
$$\a(g,x)=\a(x,g)\ \ {\mbox for\ \ all}\ \ x\in C_G(g).$$
\end{de}
Thus $g$ is $\a$-regular if and only if $\bar g \bar x=\bar x\bar g$ for all $x\in C_G(g).$

It is clear that identity element of $G$ is $\a$-regular for any cocycle $\a$.
\begin{lem}{[K]}
For any $\a\in Z^2(G,F^*)$, the following properties hold:

i) An element $g\in G$ is $\a$-regular if and only if it is $\b$-regular for any cocycle $\b$ cohomologous to $\a$.

ii) If $g\in G$ is $\a$-regular, then so is any conjugate of $g$.
\end{lem}

\begin{de}
Let $C$ be a conjugacy class of $G$ and let $g\in C$. We say that $C$ is $\a$-{\em regular} if $g$ is $\a$-regular.
\end{de}

\begin{rem}
 If we replace $\bar g$ by ${\tilde g}=\l(g)\bar g,\l(g)\in F^*,\l(1)=1$, yields 
$$ {\tilde x}{\tilde y}=\a(x,y)\l(x)\l(y)\l(xy)^{-1}\widetilde{xy}\ \ {\mbox for\ \ all }\ \ x,y\in G.$$

Hence, performing a diagonal change of basis $\{\bar g| g\in G\}$ results in replacing $\a$ by a cohomologous cocycle. The advantage of this observation is that in choosing a distinguished cocycle cohomologous to $\a$, we can work exclusively within the algebra $F^{\a}[G]$.
\end{rem}
 
\begin{de}
We say that $\a\in Z^2(G,F^*)$ is a {\em normal cocycle} if $\a(x,g)=\a(xgx^{-1},x)$ for all $x\in G$ and all $\a$-regular $g\in G$.
\end{de}

Therefore $\a$ is normal if and only if $\bar x\bar g\bar x^{-1}=\overline{xgx^{-1}}$ for all $x\in G$ and all $\a$-regular $g\in G$.

\begin{lem}\label{a2.1}{[K]}
Let $\{g_1,...,g_r\}$ be a full set of representatives for the $\a$-regular conjugacy classes of $G$ and, for each $i\in \{1,...,r\}$, let $T_i$ be a left transversal for $C_G(g_i)$ in $G$. Put $\widetilde{tg_it^{-1}}=\bar t \overline{g_i}\overline{t}^{-1},1\leq i \leq r$, and $\tilde{g}=\bar g $ if $g$ is not $\a$-regular . Then $\tilde {x}\tilde {g}\tilde {x}^{-1}=\widetilde{xgx^{-1}}$ for all $x\in G$ and all $\a$-regular $g\in G$. Thus any cocycle $\a\in Z^2(G,F^*)$ is cohomologous to a normal cocycle.
\end{lem}

\begin{th}\label{a2.2}{[K]}
Let $F$ be an arbitrary field, let $\a\in Z^2(G,F^*)$ and let $\{g_1,...,g_r\}$ be a full set of representatives for the $\a$-regular conjugacy classes of $G$. Denote by $T_i$ a left transversal for $C_G(g_i)$ in $G$ and by $C_i$ the $\a$-regular conjugacy class of $G$ containing $g_i,1\leq i\leq r$. Then the elements 
$$z_i=\sum_{t\in T_i}\bar t \overline{g_i}\overline{t}^{-1}\ \ \ (1\leq i\leq r)$$
constitute an $F$-basis $Z(F^{\a}[G])$. In particular, if $\a$ is a normal cocycle, then the elements 
$$z_i=\sum_{g\in C_i}\bar g\ \ \ (1\leq i\leq r)$$  
constitute an $F$-basis $Z(F^{\a}[G])$.
\end{th}

\section{The Generalized twisted double}\label{s3}

In this section we construct a finite dimensional semisimple
associative algebra $A_{\alpha}(G,S)$ over $\C$ associated 
to a finite group $G,$ a finite right $G$-set $\S$  and 
a suitable 2-cocycle $\alpha$. Although this construction in the present
form is influenced by the work in
[DPR], [M] and [DM5] and Section 6 of this paper, its origin 
goes back to the earlier work of [S], [BCM] and [DT] in Hopf algebra,
where it is called crossed product. 
In the case the $G$-set $\S$ consists 
of a dual basis of $\C[G]^*,$ $A_{\alpha}(G,\S)$ is exactly
the twisted double studied in [M] and [DM5]. This 
should explain why we call $A_{\alpha}(G,\S)$ a generalized 
twisted double. The exposition of this section
follows [M] closely and ideas are also similar. 

Let $G$ be any finite group with identity  $1$
and $\S$ be a finite right $G$-set. Set
$$\C{\S}=\bigoplus_{ s\in {\S}}\C e(s)$$
which is an associative algebra under product $e(s)e(t)=\d_{s,t}e(t)$
and isomorphic to $\C^{|S|}.$ The right action of $G$ on $\S$
induces a right $G$-module structure on $\C{\S}.$
Let 
$${\U}(\C{\S})=\{\sum_{ s\in {\S} }\l_{s} e(s)|\l_{s} \in \C^{*}\}$$ be the
set of units of $\C{\S}.$  It is easy to check that $\U(\C{\S})$ is a 
multiplicative right $G$-module.

From now on we will fix  an element $\a\in Z^2(G,\U(\C\S))$ such that $\a(h,1)=\a(1,h) =\sum_{s\in {\S}}e(s)$ for all $h$ in $G$. Then $\a$ defines 
 functions $\a_s: G\times G\to \C^*$ for $s\in {\S}$ such that 
$$\a(h,k)=\sum_{s\in S}\a_s(h,k)e(s)$$
for $h,k\in G.$ The cocycle condition $\a(hk,l)\a(h,k)^l=\a(h,kl)\a(k,l)$
implies immediately that  
$$\a_s(hk,l)\a_{(s)l^{-1}}(h,k)=\a_s(h,kl)\a_s(k,l)$$
for $h,k,l\in G$ and $s\in\S.$ 

The main object that we  study in this section is 
the vector space ${\A}_{\a}(G,{\S})=\C[G]\otimes \C\S$ with a basis
$g\o e(s)$ for $g\in G$, and $s\in {\S}$.

\begin{prop}

i) The ${\A}_{\a}(G,{\S})$ is an associative algebra 
under multiplication defined by 
$$ g\o e(s)\cdot h\o e(t)=\a_t(g,h)gh\o e(sh)e(t).$$

ii) If $\a$ is cohomologous to $\b$ then 
${\A}_{\a}(G,{\S})\cong{\A}_{\b}(G,{\S})$ as algebras. 
\end{prop}
\pf For i), we first prove that the product is associative. 
Let $g,h,k\in G$ and $s,t,u\in \S.$ Then 
\begin{eqnarray*}
& &\ \ \ \ (g\o e(s)\cdot h\o e(t))\cdot k\o e(u)\\
& &=(\a_t(g,h)gh\o e(sh)e(t))\cdot k\o e(u)\\
& &=\d_{sh,t}\a_t(g,h)gh\o e(t)\cdot k\o e(u)\\
& &=\d_{sh,t}\a_t(g,h)\a_u(gh,k)ghk\o e(tk)e(u)\\
& &=\d_{sh,t}\d_{tk,u}\a_t(g,h)\a_u(gh,k)ghk\o e(u)\\
& &=\d_{sh,t}\d_{tk,u}\a_{uk^{-1}}(g,h)\a_u(gh,k)ghk\o e(u)\\
& &=\d_{sh,t}\d_{tk,u}\a_{u}(g,hk)\a_u(h,k)ghk\o e(u)
\end{eqnarray*}
and
\begin{eqnarray*}
& &\ \ \ \ g\o e(s)\cdot(h\o e(t)\cdot k\o e(u))\\
& &=g\o e(s)\cdot(\a_u(h,k)hk\o e(tk)e(u))\\
& &=g\o e(s)\cdot(\d_{tk,u}\a_u(h,k)hk\o e(u))\\
& &=\d_{tk,u}\a_u(g,hk)\a_u(h,k)ghk\o e(s(hk))e(u)\\
& &=\d_{tk,u}\d_{s(hk),u}\a_u(g,hk)\a_u(h,k)ghk\o e(u).
\end{eqnarray*}
Thus $(g\o e(s)\cdot h\o e(t))\cdot k\o e(u)=g\o e(s)\cdot(h\o e(t)\cdot k\o e(u))$ and associativity holds.

A straightforward verification
\begin{eqnarray*}
& &\ \ \ \ g\o e(s)\cdot\sum_{t\in {\S}}1\o e(t)\\
& &=\sum_{t\in {\S}}\a_t(g,1)g\o e(s)e(t)\\
& &=g\o e(s)\\ 
\end{eqnarray*}
and
\begin{eqnarray*}
& &\ \ \ \ (\sum_{t\in {\S}}1\o e(t))\cdot g\o e(s)\\
& &=\sum_{t\in {\S}}\a_s(1,g)g\o e(tg)e(s)\\
& &=g\o e(s)
\end{eqnarray*}
shows that  $\sum_{t\in {\S}}1\o e(t)$ is an identity on ${\A}_{\a}(G,{\S})$.

For ii) since $\b$ is cohomologous to $\a$, 
so there exists a map $\l:G\rightarrow{\U}(\C{\S})$ 
such that $\b(x,y)=\a(x,y)\l(x)^y\l(xy)^{-1}\l(y)$ 
for all $x,y\in G$. For each $x\in G$ , we can 
rewrite $\l(x)$ as $\sum_{s\in {\S}}\l_s(x)e(s)$. 
Therefore, we have 
$\b_s(x,y)=\a_s(x,y)\l_{sy^{-1}}(x)\l_s(xy)^{-1}\l_{s}(y)$,
for any $s\in {\S}, x,y\in G.$ 
Let $f:{\A}_{\a}(G,{\S})\rightarrow{\A}_{\b}(G,{\S})$ 
defined by $f(g\o e(s))=\l_{s}(g)^{-1}(g\o e(s)).$
Since 
\begin{eqnarray*}
& &\ \ f(g\o e(s)\cdot h\o e(t))\\
& &=f(\a_t(g,h)gh\o e(sh)e(t))\\
& &=\d_{sh,t}\l_{t}(gh)^{-1}\a_t(g,h)gh\o e(t)
\end{eqnarray*}
and
\begin{eqnarray*}
& &\ \ f(g\o e(s))\cdot f(h\o e(t))\\
& &=\l_{s}(g)^{-1}g\o e(s)\cdot \l_{t}(h)^{-1}h\o e(t)\\
& &=\l_{s}(g)^{-1}\l_{t}(h)^{-1}\b_t(g,h)gh\o e(sh)e(t),
\end{eqnarray*}
${\A}_{\a}(G,{\S})$ and ${\A}_{\b}(G,{\S})$ 
are isomorphic as algebras.\qed

\begin{rem}\label{ra} {\rm The algebra $A_{\alpha}(G,\S)$ is essentially the
crossed product $\C{\S}\sharp \C[G]$ except that $\C\S$ in our setting
is a right $\C[G]$-module instead of left $\C[G]$-module in the setting
of crossed product. So one could
have used the results from [S], [BCM] and [DT] to conclude
that $A_{\alpha}(G,\S)$ is an associative algebra 
after a careful identifying two different settings. 
Since the proof
is not very long we chose to give a complete proof.} 
\end{rem}

\begin{rem} {\em If we take ${\S}=\{e(g)|g\in G\}$ to be the
dual basis of $\C[G]^*$ which is a right $G$-set via
$e(g)\cdot h=e(h^{-1}gh)$ for $g,h\in G$ we get the associative
algebra $D_{\alpha}(G)$ in [DM5] which is a deformation
of the Drinfeld's double $D(G).$ It turns out
that there exists $\alpha$ determined naturally 
by the ``twisted representations'' such that $D_{\alpha}(G)$ is the right 
algebra in the study of holomorphic orbifold conformal field theory
(see [DVVV] and [DM5]).} 
\end{rem} 
 
For each $s\in {\S}$, let $G_s=\{h\in G|s\cdot h=s\}$ be the stabilizer
of $s.$ Note that $Res_{G_s}^{G}\a_s$ is in $Z^2(G_s,\C^*).$
Let ${\O}_s$ be the orbit of $s$ under $G$ and 
$G=\bigcup_{i=1}^kG_sg_i$ be a right coset decomposition with $g_1=1.$ 
Then 
${\O}_s=\{sg_i|i=1,...,k\}$ and $G_{s\cdot g_i}=g_i^{-1}G_sg_i.$ We
define several subspaces of  ${\A}_{\a}(G,{\S}):$
$$S(s)=\<a\o e(s)|a\in G_{s}\>\ \ ,\ \ N(s)=\<a\o e(s)|a\in G\setminus G_{s}\>,$$
$$
D(s)=\<a\o e(s)|a\in G\>\ \ ,\ \ 
D({\O}_s)=\<a\o e(s\cdot g_i)|i=1,...,k, a\in G\>.$$ 
Then $D(s)=S(s)\oplus N(s)$.

Decompose $\S$ into a disjoint union of orbits ${\S}=\bigcup_{j\in J}{\O}_j$.
Let $s_j$ be a representative element of ${\O}_j.$ Then
${\O}_j=\{s_j\cdot h|h\in G\}$ and $ {\A}_{\a}(G,{\S})=\bigoplus_{j\in J}D({\O}_{s_j}).$ 

\begin{lem}\label{l3.2} Let $s\in \S$ and $G=\cup_{i=1}^kG_sg_i.$ Then

1) $S(s)$ is a subalgebra of ${\A}_{\a}(G,{\S})$ 
isomorphic to $\C^{\a_{s}}[G_s]$ where $\C^{\a_{s}}[G_s]$ is the
twisted group algebra.

2) $N(s)$ is a 2-sided nilpotent ideal of $D(s)$ and $D(s)\cdot
N(s)=0.$ 

3) $D({\O}_s)=\bigoplus_{i=1}^kD(sg_i)$ is a direct sum of left ideals.

4) Each $D({\O}_s)$ is a 2-sided ideal of ${\A}_{\a}(G,{\S})$
and ${\A}_{\a}(G,{\S})=\oplus_{j\in J}D({\O}_{s_j}).$ Moreover,
$D({\O}_s)$ has identity element $\sum_{t\in {\O}_s} 1\otimes e(t).$ 
\end{lem} 
\pf 
1) Clearly, $S(s)$ is a subalgebra of ${\A}_{\a}(G,{\S}).$ 
Let $\r:S(s)\rightarrow \C^{\a_{s}}[G_s]$ be a linear isomorphism
 determined
by $\r(a\o e(s))=a$. Then
\begin{eqnarray*}
& & \ \ \ \ \r(a\o e(s)\cdot b\o e(s))\\
& &=\r(\a_{s}(a,b)ab\o e(sb)e(s))\\
& &=\r(\a_{s}(a,b)ab\o e(s))\\
& &=\a_{s}(a,b)ab\\
& &=\r(a\o e(s))\r(b\o e(s)).
\end{eqnarray*}
So, $\r$ is an algebra isomorphism.

2) Since for $b\in G\setminus  G_{s},$ $sb\ne s$. We immediately
have $ a\o e(s)\cdot b\o e(s)=\a_{s}(a,b)ab\o e((s)b)e(s)=0$
for $a\o e(s)\in D(s)$ and $b\o e(s)\in N(s)$. 

It remains to show  that $N(s)$ is right ideal. 
Take $a\o e(s)\in N(s)$ and $b\o e(s)\in S(s).$  
Then $a\o e(s)\cdot b\o e(s)=\a_{s}(a,b)ab\o e((s)b)e(s)=\a_{s}(a,b)ab\o e(s).$ If $ab$ lies in $G_{s}$ so does $a$. Thus $ab\in G\setminus G_{s}$ and $ a\o e(s)\cdot b\o e(s)\in N(s)$.

3)-4) are clear.  \qed 

Let $A$ be an algebra (eg., associative algebra, Lie algebra or vertex
operator algebra). We denote the module category of $A$ by $A$-Mod. 

For convenience we set $d_{g,s}=g\otimes e(s)$ for $g\in G$ and $s\in S.$

\begin{th}\label{t3.3} The functors
\begin{eqnarray*}
 M&\stackrel{f}{\mapsto}D(s)&\bigotimes_{S(s)}M\\
N&\stackrel{g}{\mapsto}&d_{1,s}N\\
\end{eqnarray*}
($M\in \C^{\a_s}[G_s]$-{\em Mod}, $N\in D({\O}_s)$-{\em Mod})
define an equivalence between categories $\C^{\a_s}[G_s]$-{\em Mod} and $D({\O}_s)$-{\em Mod}. In particular, the simple $\C^{\a_s}[G_s]$-modules are mapped to simple $D({\O}_s)$-modules and conversely. 
\end{th}
\pf Recall that $G=\bigcup_{i=1}^k G_s g_i$ be a coset decomposition with $g_1=1.$  Then 
$G=\bigcup_{i=1}^k g_i^{-1}G_s.$
Note that $g_i^{-1}a\otimes e(s)$ can be 
rewritten as 
$\alpha_s(g_i^{-1},a)^{-1}g_i^{-1}\otimes e(s)\cdot a\otimes e(s)$
for any $i$ and $a\in G_s.$  
 Let $M\in \C^{\a_s}[G_s]$-Mod. Then an arbitrary vector in $D(s)\bigotimes_{S(s)}M$ has expression 
$\sum_{i=1}^k d_{g_i^{-1},s}\o_{S(s)} m_i$ for some $m_i\in M.$ 
From 
$$d_{1,s}(\sum_{i=1}^kd_{g_i^{-1},s}\o_{S(s)} m_i)=d_{1,s}\o_{S(s)} m_1$$
we see that $d_{1,s}D(s)\bigotimes_{S(s)}M$ is a subset of $d_{1,s}\otimes_{S(s)}M.$ On the other hand, $d_{1,s}\o_{S(s)} m= d_{1,s}(d_{1,s}\o_{S(s)} m)$
 belongs to $d_{1,s}(D(s)\bigotimes_{S(s)}M)$. Thus $d_{1,s}\otimes_{S(s)}M$
and $d_{1,s}(D(s)\bigotimes_{S(s)}M)$ are equal.  Since $d_{1,s}\bigotimes_{S(s)}M$ is isomorphic to $M$ as
$\C^{\a_s}[G_s]$-modules . This implies that $g\circ f(M)$ and $M$ are
isomorphic as $\C^{\a_s}[G_s]$-modules.

Now take $N\in D({\O}_s)$-Mod. We are going to show that $D(s)\otimes_{S(s)}d_{1,s}N$ and $N$ are isomorphic as $D({\O}_s)$-modules. For this purpose we define a linear map
\begin{eqnarray*}
\r: D(s)\bigotimes_{S(s)} d_{1,s}N &\rightarrow & N\\
  d_{b,s}\o_{S(s)}d_{1,s}n&\mapsto& d_{b,s}n.
\end{eqnarray*}
Noting that $\r(d_{b,s}\o_{S(s)}d_{1,s}n)=d_{b,s}d_{1,s}n$ we immediately
see that $\r$ is well defined. 
The fact that $\r$ is a $D({\O}_s)$-module homomorphism follows from
\begin{eqnarray*} 
& &\ \ \ \ \r(d_{a,sg_i}d_{b,s}\o_{S(s)}d_{1,s}n)\\
& &=\r(\d_{sg_ib,s}\a_s(a,b)d_{ab,s}\o_{S(s)} d_{1,s}n)\\
& &=\d_{sg_ib,s}\a_s(a,b)d_{ab,s}n\\
\end{eqnarray*}
and
\begin{eqnarray*}
& &\ \ \ \ d_{a,sg_i}\r (d_{b,s}\o_{S(s)}d_{1,s}n)\\ 
& &=d_{a,sg_i}d_{b,s}n\\
& &=\d_{sg_ib,s}\a_s(a,b)d_{ab,s}n
\end{eqnarray*}
for $a,b\in G,$ $1\leq i\leq k$ and $n\in N.$ 

It remains to show that $\r$ is a bijection.  As $\sum_{i=1}^kd_{1,sg_i}$
is the identity of $D(\O_s),$ we have
$n=\sum_{i=1}^k d_{1,sg_i}n$.
For fixed $i$ one can easily see that 
$$\r(d_{g_i^{-1},s}\o_{S(s)}\a_{sg_i}(g_i^{-1},g_i)^{-1}d_{g_i,sg_i}n)=d_{1,sg_i}n.$$ 
Thus $\r$ is onto.

Since $d_{1,sg_i}d_{1,sg_j}=\delta_{i,j}d_{1,sg_i}$ we have 
$$N=\bigoplus_{i=1}^kd_{1,sg_i}N.$$ Let $b=g_i^{-1}a\in g_i^{-1}G_s.$
Then $\r(d_{b,s}\o_{S(s)}d_{1,s}n)= d_{1,sg_i}d_{b,s}n\in d_{1,sg_i}N$
for $n\in N.$ In order to prove that $\r$ is one to one, it is enough
to show that if $\r(\sum_{p}d_{b_p,s}\o_{S(s)}d_{1,s}n_p)=0$ for some
$b_p\in g_i^{-1}G_s$ and $n_p\in N,$ then
$\sum_{p}d_{b_p,s}\o_{S(s)}d_{1,s}n_p=0.$ Using the relation
$d_{g_i^{-1}a,s}=\a_s(g_i^{-1},a)^{-1}d_{g_i^{-1},s}d_{a,s}$ for $a\in G_s$
we can rewrite $\sum_{p}d_{b_p,s}\o_{S(s)}d_{1,s}n_p$ as
$d_{g_i^{-1},s}\o_{S(s)}d_{1,s}n$ for some $n\in N.$ Thus
$\r(d_{g_i^{-1},s}\o_{S(s)}d_{1,s}n)=d_{g_i^{-1},s}n=0.$ Applying
$\a_s(g_i,g_i^{-1})^{-1}d_{g_i,sg_i}$ to $d_{g_i^{-1},s}n$ yields
$d_{1,s}n=0.$ This shows that $\r$ is one to one. This completes the
proof that $\r$ is an isomorphism.  \qed

\begin{th}\label{t50} We have

i)  Algebra $D({\O}_s)$ is semisimple for $s\in \S$ and simple
$D({\O}_s)$-modules are  precisely $\Ind_{S(s)}^{D(s)}M$ where $M$ ranges over the simple $\C^{\a_{s}}[G_{s}]$-modules.

ii) $A_{\a}(G,{\S})$ is semisimple and  simple $A_{\a}(G,{\S})$-modules are precisely $\Ind_{S(s_j)}^{D(s_j)}M$ where $M$ ranges over the simple $\C^{\a_{s_j}}[G_{s_j}]$-modules and $j\in J.$ 
\end{th}

\pf i) follows from Theorem \ref{t3.3} and that fact that  
$\C^{\a _s}[G_s]$-Mod is a semisimple category. ii) follows from
Lemma \ref{l3.2} and i).
\qed 

\begin{rem}{If we only wanted to 
know the semisimplicity of $A_{\alpha}(G,\S),$ we could 
find it in the literature of Hopf algebra 
when we regard $A_{\alpha}(G,\S)$ as a crossed product (see Remark \ref{ra}).
But in the later sections
we need to know the explicit structure of simple $A_{\alpha}(G,\S)$-modules
given in Theorem \ref{t3.3}. The semisimplicity of $A_{\alpha}(G,\S)$
is a trivial corollary of Theorem \ref{t3.3}.}
\end{rem}

Next we determine the center $Z(D({\O}_s))$ of $D({\O}_s)$. 
However, we are not going to use it in this paper 
but still include it for general interest. We certainly expect that
the result on center will be used in our future study on general
orbifold conformal field theory.

Since any cocycle in $Z^2(G_s,\C^*)$ is cohomologous to a normal cocycle (cf. Lemma \ref{a2.1}). Then we may {\em assume $\a_s$ is a normal cocycle} for
all $s\in S.$  We also assume that $\a_{sg_i}(g_j^{-1}hg_i,g_i^{-1}ag_i)=\a_{sg_i}(g_j^{-1}hah^{-1}g_j,g_j^{-1}hg_i)$ for  all $1\leq i,j\leq k,$   $h\in 
G_s$ and $\a_s$-regular $a\in G_s$. In the orbifold theory, these
conditions are satisfied by chosing $\alpha_s$ carefully.

Let $\{l_1,...,l_r\}$ be a full set of representatives for the $\a_s$-regular conjugacy classes of $G_s$ and for each $t$ in $\{1,...,r\}$, let $L_t$ be $\a_s$-regular conjugacy class of $G_s$ containing $l_t,1\leq t\leq r.$
Set
$$Z(L_t)=\sum_{a\in L_t}\sum_{i=1}^k g_i^{-1}ag_i\o e(sg_i).$$
\begin{lem} $Z(L_t)$ is an center element of $D(\O_s).$
\end{lem}

\pf Let $b\o e(sg_j)\in D({\O}_s)$. Then there exists $1\leq i'\leq k$ and
$h\in G_s$ such that $b=g_{i'}^{-1}hg_j.$ 
We have
\begin{eqnarray*}
& &\ \ \ \ b\o e(sg_j)\cdot Z(L_t)\\
& &=b\o e(sg_j)\cdot \sum_{a\in L_t}\sum_{i=1}^k g_i^{-1}ag_i\o e(sg_i)\\
& &=\sum_{a\in L_t}\sum_{i=1}^k\a_{sg_i}(b,g_i^{-1}ag_i)bg_i^{-1}ag_i\o e(sg_jg_i^{-1}ag_i)e(sg_i)\\
& &=\sum_{a\in L_t}\a_{sg_j}(b,g_j^{-1}ag_j)bg_j^{-1}ag_j\o e(sg_j)\\
& &=\sum_{a\in L_t}\a_{sg_j}(g_{i'}^{-1}hg_j,g_j^{-1}ag_j)g_{i'}^{-1}hag_j\o e(sg_j)
\end{eqnarray*}
and 
\begin{eqnarray*}
& &Z(L_t)\cdot b\o e(sg_j)\\
& &=\sum_{a\in L_t}\sum_{i=1}^k g_i^{-1}ag_i\o e(sg_i)\cdot b\o e(sg_j)\\
& &=\sum_{a\in L_t}\sum_{i=1}^k\a_{sg_j}(g_i^{-1}ag_i,b)g_i^{-1}ag_ib\o e(sg_ib)e(sg_j)\\
& &=\sum_{a\in L_t}\a_{sg_j}(g_{i'}^{-1}ag_{i'}, g_{i'}^{-1}hg_j)g_{i'}^{-1}ahg_j\o e(sg_j).
\end{eqnarray*}
Now the assertion follows from 
\begin{eqnarray*}
& &\ \ \ \ \sum_{a\in L_t}\a_{sg_j}(g_{i}^{-1}hg_j,g_j^{-1}ag_j)hah^{-1}\\
& &=\sum_{a\in L_t}\a_{sg_j}(g_{i}^{-1}hah^{-1}g_{i},g_{i}^{-1}hg_j)hah^{-1}\\
& &=\sum_{a\in L_t}\a_{sg_j}(g_{i}^{-1}ag_{i},g_{i}^{-1}hg_j)a.
\end{eqnarray*}
\qed

\begin{th}
Let $\{l_1,..., l_r\}$ be a full set of representatives for the $\a_s$-regular conjugacy classes of $G_s$ and for each $t\in \{1,...,r\}$, let $L_t$ be a $\a_s$-regular conjugacy class of $G_s$ containing $l_t, 1\leq t\leq r$. 
Then the elements $Z(L_t)$ constitute a $\C$-basis of $Z(D({\O}_s))$.
\end{th}
\pf By Theorem \ref{t3.3} the dimension of $Z(D({\O}_s))$ equals to
the number of inequivalent irreducible $\C^{\a_s}[G_s]$-modules which 
is $r$ (see Theorem \ref{a2.2}).
\qed

The following corollary is immediate.
\begin{cor} Let $\{l_1^{j},...,l_{r_j}^{j}\}$ be a full set of representatives for the $\a_{s_j}$-regular conjugacy classes of $G_{s_j}$ and for each  $t_j\in \{1,...,r_j\}$, let $L_{t_j}$ be $\a_{s_j}$-regular conjugacy class of $G_{s_j}$ containing $l_{t_j}, 1\leq t_j\leq r_j.$ Set
$$Z(L_{t_j})=\sum_{a\in L_{t_j}}\sum_{i_j=1}^{k_j}g_{i_j}^{-1}ag_{i_j}\o e(sg_{i_j}).$$
Then for all $j\in J$, for all $t_j\in \{ 1,..., r_j\}$, $Z(L_{t_j})$ constitutes a $\C$-basis for $Z(A_{\a}(G,{\S})).$
\end{cor}

\section{Modules for vertex operator algebras and related results}

In this section we turn our attention to the theory of vertex operator 
algebras. In particular we shall defines various notion of modules
for a vertex operator algebra $V$ following [FLM], [DLM2] and [Z]. We also
recall from [DLM4] the associative algebras $A_n(V)$ for any
nonnegative integer $n$ and relevant results. These results will be used
extensively in Sections 5 and 6 to study dual pairs arising from
an action of a finite group $G$ on $V.$  

Let $V=(V,Y,{\bf 1},\omega)$  be a vertex operator algebra
(see [B] and [FLM]).  
For a vector space $W$, let $W\{z\}$ be the space of $W$-valued formal series in arbitrary complex powers of $z$. We present three different
notion of modules (cf. [FLM], [DLM2] and [Z]). 
\begin{de} A {\em weak $V$-module} $M$ is a vector space equipped 
with a linear map
$$\begin{array}{lll}
V&\to &(\End\,M)\{z\}\\
v&\mapsto&\displaystyle{ Y_M(v,z)=\sum_{n\in\Q}v_nz^{-n-1}\ \ \ (v_n\in
\End\,M)}
\end{array}$$
which satisfies the following properties for all $u\in V$, $v\in V,$ 
$w\in M$,
\begin{eqnarray}
& &u_lw=0\ \ \  				
\mbox{for}\ \ \ l>>0\label{vlw0}\\
& &Y_M({\mathbf 1},z)=1;\label{vacuum}
\end{eqnarray}
 \begin{equation}\label{jacobi}
\begin{array}{c}
\displaystyle{z^{-1}_0\delta\left(\frac{z_1-z_2}{z_0}\right)
Y_M(u,z_1)Y_M(v,z_2)-z^{-1}_0\delta\left(\frac{z_2-z_1}{-z_0}\right)
Y_M(v,z_2)Y_M(u,z_1)}\\
\displaystyle{=z_2^{-1}
\delta\left(\frac{z_1-z_0}{z_2}\right)
Y_M(Y(u,z_0)v,z_2)},
\end{array}
\end{equation}
where $\delta(z)=\sum_{n\in\Z}z^n$ and all binomial expressions (here and below) are to be expanded in nonnegative integral powers of the second variable. Elementary properties of the $\delta$-function can be found in [FLM] and [FHL]
\end{de}

(\ref{jacobi}) is called the {\em Jacobi identity}. One can prove 
that the Jacobi identity is equivalent to the following
associativity formula (see [FLM] and [FHL])
\begin{eqnarray}\label{ea}
(z_{0}+z_{2})^{k}Y_{M}(u,z_{0}+z_{2})Y_{M}(v,z_{2})w
=(z_{2}+z_{0})^{k}Y_M(Y(u,z_0)v,z_2)w
\end{eqnarray}
and commutator relation
\begin{eqnarray}
& &\ \ \ \  [Y_{M}(u,z_{1}),Y_{M}(v,z_{2})]\nonumber\\
& &=\Res_{z_{0}}z_2^{-1}
\delta\left(\frac{z_1-z_0}{z_2}\right)Y_M(Y(u,z_0)v,z_2)\label{ec}
\end{eqnarray}
where $w\in M$ and $k$ is a nonnegative integer such that $z^{k}Y_{M}(u,z)w$ involves only nonnegative integral powers of $z.$

We may also deduce from
(\ref{vlw0})-(\ref{jacobi}) the usual Virasoro algebra axioms (see [DLM2]).
Namely, 
if $Y_M(\om,z)=\sum_{n\in\Z}L(n)z^{-n-2}$ then
\begin{equation}
\label{g3.8}
[L(m),L(n)]=(m-n)L(m+n)+\frac{1}{12}(m^3-m)\delta_{m+n,0}(\mbox{rank}\,V)
\end{equation}
and 
\begin{equation}\label{2.5}
\frac{d}{dz}Y_M(v,z)=Y_M(L(-1)v,z).
\end{equation}

Suppose that $(M_{i},Y_{i})$ are two weak $V$-modules, $i$=1,2.
A homomorphism from $M_{1}$ to $M_{2}$ is a linear map
$f$:$M_{1}\rightarrow M_{2}$ which satisfied
$fY_{M_{1}}(v,z)=Y_{M_{2}}(v,z)f$ for all $v\in V$. 
We call $f$ an isomorphism if $f$ is also a linear isomorphism.

\begin{de}\label{d1} A (ordinary) {\em $V$-module} is a weak $V$-module $M$ which carries a $\C$-grading induced by the spectrum of $L(0)$. Then 
$M=\bigoplus_{\lambda\in \C}M_{\lambda}$ where 
$M_{\lambda}=\{w\in M|L(0)w=\lambda w\},$ dim$M_{\lambda}<\infty.$
Moreover, for fixed $\lambda, M_{n+{\lambda}}=0 $ for all small enough integers $n.$
\end{de}

The notion of module here is essentially the notion of module given
in [FLM].

\begin{de}\label{d2} An {\em admissible $V$-module} is a weak 
$V$-module $M$ which carries a $\Z_{+}$ grading
$M=\bigoplus_{n\in \Z_{+}}M(n)$ satisfying the following
condition:$$v_{m}M(n)\subset M(n+\wt v-m-1)$$ for homogeneous $v\in
V,$ where $\Z_+$ is the set of the nonnegative integers. 
\end{de}

The notion of admissible module here is the notion of module 
in [Z]. Using a grading shift we can always arrange the 
grading on $M$ so that $M(0)\ne 0.$ This shift is important
in the study of algebra $A_n(V)$ below.

It is not too hard to see that any $V$-module is an admissible $V$-module. 
So there is a natural identification of the category
of $V$-modules with a subcategory of the category
of admissible $V$-modules.

\begin{de}\label{ration} $V$ is called rational if every admissible $V$-module
is a direct sum of irreducible admissible $V$-modules.
\end{de}

It is proved in [DLM3] that if $V$ is rational then there are only
finitely many inequivalent irreducible admissible $V$-modules
and each irreducible admissible $V$-module is an ordinary module.

The following proposition can be found  in [L] and [DM2].
\begin{prop}\label{4.19} If $M$ is a simple weak $V-$module then
$M$ is spanned by $\{u_nm|u\in V,n\in \Q\}$ where $m\in M$ is a fixed nonzero 
vector.
\end {prop}

We now recall the associative algebra $A_{n}(V)$ as constructed in [DLM4]. 
Let $O_n(V)$ be the linear span of all $u\circ_n v$ and $L(-1)u+L(0)u$
where for homogeneous $u\in V$ and $v\in V,$
\begin{equation}\label{g2.2}
u\circ_n v=\Res_{z}Y(u,z)v\frac{(1+z)^{\wt u+n}}{z^{2n+2}}.
\end{equation}
Define the linear space $A_n(V)$ to be the quotient $V/O_{n}(V).$ 
We also define a second product $*_n$ on $V$ for $u$ and $v$ as
above: 
\begin{equation}\label{a5.1}
u*_nv=\sum_{m=0}^{n}(-1)^m{m+n\choose n}\Res_zY(u,z)\frac{(1+z)^{\wt\,u+n}}{z^{n+m+1}}v.
\end{equation} 
Extend linearly to obtain a bilinear product  on $V$.

Let $M=\sum_{n\in\Z_{+}}M(n)$ be an admissible $V$-module. Following
[Z] we define weight zero operator $o_M(v)=v_{\wt v-1}$ on $M$ for
homogeneous $v$ and extend $o_M(v)$ to all $v$ by linearity. It is
clear from the definition that $o_M(v)M(n)\subset M(n)$ for
all $n.$  

\begin{th}\label{l2.3}  Let $V$ be a vertex operator algebra and
$M$ an admissible $V$-module. Then

1) The product $*_{n}$ induces the structure of an 
associative algebra  on $A_{n}(V)$ with identity ${\bf 1}+O_{n}(V).$ 
Moreover $\omega+O_{n}(V)$ is a central element of $A_{n}(V).$

2) For $0\leq m\leq n,$ the map $\psi_n: v+O_n(V)\mapsto o_M(v)$ from
$A_n(V)$ to $\End M(m)$ makes $M(m)$ an $A_n(V)$-module.

3) $M$ is irreducible if and only if $M(n)$ is a simple $A_n(V)$-module
for all $n.$

4) The identity map on $V$ induces an onto algebra homomorphism
from $A_n(V)$ to $A_m(V)$ for $0\leq m\leq n.$ 

5) Two irreducible admissible $V$-modules $M^1$ and $M^2$ 
with $M^1(0)\ne 0$ and $M^2(0)\ne 0$ are isomorphic if and only if
$M^1(n)$ and $M^2(n)$ are isomorphic $A_n(V)$-modules 
for any $n\geq 0$.
\end{th}

Parts 1)-2) and 4)-5) are proved in [DLM4], and 3) is given 
in [DM4]. 

\section{Dual pair I}

In the next two sections we assume that $V$ is a simple vertex operator
algebra. Our main goal is to generalize the duality result obtained
in [DLM1] on $V$ to an arbitrary irreducible $V$-module. More precisely,
let $G$ be a finite automorphism group of $V$ and $M$ an irreducible
$V$-module. Then $M$ is a module
for the vertex operator subalgebra $V^G$ of $G$-invariants. There are three
 questions
one can ask: (1) Is $M$ is a completely reducible $V^G$-module?
(2) What are irreducible modules for $V^G$ which occur as submodules
of $M?$ (3) What are the relation between irreducible $V^G$-submodules
of $M$ and $V^G$-submodules of another irreducible $V$-module $N?$
These questions will be answered completely in this section
and the next section.

In this section we deal with the case that $M$ is $G$-stable (see the
definition below). The general case will be treated in the next section.
The basis ideas in the proof of main theorems come from [DLM1].

Let $(M,Y_M)$ be an irreducible $V$-module 
and $g\in G$. Following [DM1] we set $M\circ g=M$ as vector spaces, 
and $Y_{M\circ g}(v,z)=Y_M(gv,z)$. 
Note that $M\circ g$ 
is also an irreducible $V$-module. We {\em assume} in this section 
that $M$ is $G$-{\em stable} in the sense that for any $g\in G,$
$M\circ g$ and $M$ are isomorphic.  
If $M=V$ this assumption is always true. 

For $g\in G$ there is a linear isomorphism $\phi(g): M\to
M$ satisfying
\begin{equation}\label{3.4'}
\phi(g)Y_M(v,z)\phi(g)^{-1}=Y_{M\circ g}(v,z)=Y_M(gv,z)
\end{equation}
for $v\in V.$ The simplicity of $M$ together with
Schur's lemma shows that $g\mapsto \phi(g)$ is a projective
representation of $G$ on $M.$  
Let $\a_M$  be the corresponding 2-cocycle in $C^2(G,\C^{\times}).$
Then $M$ is a module for $\C^{\a_M}[G]$
where $\C^{\a_M}[G]$ is the twisted group algebra.
It is worth pointing out that if $M=V$ we can take $\phi(g)=g$
and $\C^{\a_M}[G]=\C[G].$

Let $\Lambda_{G,\a_M}$ be the set of all irreducible characters $\lambda$ 
of  $\C^{\a_M}[G]$. We denote
the corresponding simple module by $W_{\l}.$ Again if $M=V$,  
$\Lambda_{V,\a_V}=\hat G$ is the set of all irreducible characters of $G.$ 
Note that $M$ is a semisimple 
$\C^{\a_M}[G]$-module. 
Let $M^{\lambda}$ be the sum of simple 
$\C^{\a_M}[G]$-submodules of $M$ isomorphic 
to $W_{\l}.$ Then $M=\oplus_{\lambda\in \Lambda_{G,\a_M}}M^{\lambda}$. 
Moreover, $M^{\lambda}=W_{\lambda}\otimes M_{\l}$ 
where $M_{\l}=\hom_{\C^{\a_M}[G]}(W_{\l},M)$ is the multiplicity
of $W_{\l}$ in $M.$ 
As in [DLM1], we can, 
in fact, realize $M_{\l}$ as
a subspace of $M$ in the following way.  Let $w\in W_{\l}$ 
be a fixed nonzero vector. Then we can identify 
$\hom_{\C^{\a_M}[G]}(W_{\l},M)$ with the subspace
$$\{f(w) |f\in \hom_{\C^{\a_M}[G]}(W_{\l},M)\}$$
 of $M^{\l}.$  

We also set $V^G=\{v\in V|gv=v, g\in G\}.$ Then $V^G$ is a simple vertex
operator subalgebra of $V$ (see [DM2]) 
and $\phi(g)Y_M(v,z)\phi(g)^{-1}=Y_M(gv,z)
=Y_M(v,z)$ for all $v\in V^G.$ Thus the actions of $\C^{\a_M}[G]$
 and $V^G$ on $M$ are commutative. This shows that
both $M^{\l}$ and $M_{\l}$ are ordinary $V^G$-modules.

\begin{lem}\label{corect} If $M_{\l}\ne 0$ then $M_{\l}$ is an irreducible
$V^G$-module.
\end{lem}

\pf By Theorem \ref{l2.3} 3) it is enough to prove that 
$M_{\l}(n)=M_{\l}\cap M(n)$ is a simple $A_n(V^G)$-module for
all $n\geq 0.$ It is equivalent to show 
that for any $n\in \Z, n\geq 0$ 
$M_{\lambda}(n)$ is generated by 
any nonzero vector in $M_{\lambda}(n)$
as an $A_{n}(V^G)$-module. 

Note from the definition of $A_n(V)$ that $A_n(V)$ is a $G$-module via 
$g(v+O_n(V))=gv+O_n(V)$ for $g\in G$ and $v\in V.$  Also $\End M(n)$ is an
$G$-module via $gf=\phi(g)f\phi(g)^{-1}$ for $g\in G$ and $f\in
\End M(n).$ Then the algebra homomorphism $v+O_n(V)\mapsto o_M(v)$
from $A_n(V)$ to $\End M(n)$ is a $G$-homomorphism as 
$\phi(g)o_M(v)\phi(g)^{-1}=o_M(gv).$ From Proposition \ref{4.19} 
and Theorem \ref{l2.3} 3) we see that $\End M(n)=\{o_M(v)|v\in V\}.$
Since $G$ is a finite group we immediately have 
$\psi(A_n(V)^G)=(\End M(n))^G=\{o_M(v)|v\in V^G\}.$  
 
Let $x,y\in M_{\l}(n)$ be linearly independent. 
Since $M(n)$ is a $\C^{\a_M}[G]$-module, 
we can write $M(n)$ as a direct sum
$W_{\lambda}\otimes x\oplus W_{\lambda}
\otimes y\oplus W $ where $W$ 
is a $\C^{\a_M}[G]$-submodule of $M(n).$ 
Define a map $\beta\in\End M(n)$ 
such that $\beta(u\otimes x+v\otimes y+w)=
v\otimes x+u\otimes y+w$ for $u,v\in W_{\lambda}$ 
and $w\in W.$ Then $\beta\in(\End M(n))^{G}.$ 
Thus there exists vectors $v\in V^G$ such that
$o_M(v)x=y.$ This implies that $M_{\l}(n)$ is an simple $A_n(V)^G$-module.
Since $M_{\l}(n)$ is also an $A_n(V^G)$-module with the same action,
we see that $M_{\l}(n)$ is a simple $A_n(V^G)$-module. This can be also
explained by noting that the identity on $V^G$ induces an onto 
algebra homomorphism from $A_n(V^G)$ to $A_n(V)^G.$ 
\qed

\begin{rem}\label{rcorect}{\rm If we take $M=V$ in Lemma \ref{corect}, then
each $V_{\l}$ is an irreducible $V^G$-module for any irreducible character 
$\l$ of $G.$ This is the part (1) of Theorem 2.4 in [DLM1],
where Lemma 2.2 of [DLM1] is used. It is pointed out to us by Wassermann
that there is a gap in Lemma 2.2 of [DLM1]. The special case of Lemma
\ref{corect} now fixes the problem in [DLM1]. (The group $G$
in [DLM1] can be a compact group. But the correction present here
works for such $G$.) The proofs for parts (2)
and (3) of Theorem 2.4 in [DLM1] are correct.}
\end{rem}

For the purpose of continuous discussion we now recall 
Theorem 2.4 from [DLM1] (also see Remark \ref{rcorect}).

\begin{th}\label{t1y} Let $V$ be a simple vertex operator algebra and
$G$ a finite automorphism group. Then

1) $V^{\l}$ is nonzero for any $\l\in \hat G.$ 

2) Each $V_{\l}$ is an irreducible $V^G$-module.

3) $V_{\l}$ and $V_{\mu}$ are equivalent $V^G$-module if and only
if $\l=\mu.$
\end{th}

The main result in this  section is the following. 

\begin{th}\label{t2y} With the same notation as above we have:

1) $M^{\l}$ is nonzero for any $\l\in \Lambda_{G,\a_M}.$

2) Each $M_{\l}$ is an irreducible $V^G$-module.

3) $M_{\l}$ and $M_{\gamma}$ are equivalent $V^G$-module if and only
if $\l=\gamma.$

That is, $V^G$ and $\C^{\a_M}[G]$ form a dual pair on $M$ in the sense
of Howe (see [H1] and [H2]).
\end{th}

Let $\gamma\in \hat G$ and $\mu \in \Lambda_{G,\a_M}$
such that $M^{\mu}\ne 0.$ By Theorem \ref{t1y} 
we can regard $W_{\gamma}$ and
$W_{\mu}$ as a $\C[G]$-submodule of $V$ and a $\C^{\a_M}[G]$-submodule of $M$, respectively.  In fact, we can assume
that $W_{\gamma}$ and $W_{\mu}$ are homogeneous
subspaces of $V$ and $M,$ respectively. Following
[DM3] we define 
a subspace $Z_s$ of $M$ for any $s\in \Z$ by
$$Z_s=\{\sum_{m\geq s}u_mw|u\in W_{\gamma}, w\in W_{\mu}\}.$$ Then it
is easy to see that $Z_s$ is a $\C^{\a_M}[G]$-submodule of $M.$ We also
define a map
$$\psi_s: W_{\gamma}\otimes W_{\mu}\to Z_s$$
such that $\psi_s(u\otimes w)=\sum_{m\geq s}u_mw.$ Recall twisted 
group algebra from Section 2. In particular, $\{\bar g|g\in G\}$
is a basis of $\C^{\a_M}[G].$ Then 
$W_{\gamma}\otimes W_{\mu}$ is a 
$\C^{\a_M}[G]$-module with $\bar g$ acting as $g\otimes \phi(g)$
and  $\psi_s$ is  a $\C^{\a_M}[G]$-module homomorphism. 

\begin{lem}\label{l4.2y} For all sufficiently small $s,$ $\psi_s$
is a $\C^{\a_M}[G]$-module isomorphism.
\end{lem}

This result in the case $M=V$ is proved in [DM2] and [DM3]. The proof
is similar in the general case. Here we give an outline of the proof
and we refer the reader to [DM2] and [DM3] for the details. From the
associativity formula (\ref{ea}) and commutator formula (\ref{ec}) we
can prove the following fact. Let $u^i\in V$ for $i=1,...,n$ and
$w^1,...,w^n\in M$ be linearly independent. Then
$\sum_{i=1}^nY_M(u^i,z)w^i\ne 0$ (see the proof of Lemma 3.1 of [DM2];
also see Proposition 11.9 of [DL]). Then following the proof of Lemma
2.2 of [DM3] for the case $M=V$ gives the result.

\bigskip

We now prove Theorem \ref{t2y}. Part 2) is Lemma \ref{corect}. 
In order to prove 1) we first note that there exists 
$\mu\in\Lambda_{M,\a_M}$ such that $M^{\mu}\ne 0.$ 
Let $\mu^*$ be the character of $\C^{\a_M^{-1}}[G]$ dual to $\mu.$ 
That is, the corresponding $\C^{\a_M^{-1}}[G]$-module
is exactly $W_{\l}^*=\Hom_{\C}(W_{\l},\C).$ 
Then $\mu ^* \otimes \lambda$ is a 
character  of $G$ for any $\lambda\in \Lambda_{G,\a_M}.$ 
 Let $\gamma$ be an irreducible character of $G$
such that $\gamma\in \mu^*\otimes \lambda.$ Thus
$$\Hom_{\C^{\a_M}[G]}(\lambda ,\mu \otimes \gamma)=\Hom_{\C[G]}(\mu^*\otimes \lambda,\gamma)\neq 0.$$
By Lemma \ref{l4.2y}, for small enough $s,$ the submodule
$Z_s$ of $M$ contains a submodule isomorphic to $W_{\l}.$ That is,
$M^{\lambda}\ne 0$ for all $\l\in \Lambda_{G,\a_M}.$

Finally we prove 3). Let $\l,\gamma\in \Lambda_{G,\a_M}$ are different. 
We can take $n\in \Z$ such that $M^{\l}(n)=M^{\l}\cap M(n)\ne 0.$
Then $M(n)$ is a direct sum of  
$\C^{\a_M}[G]$-modules
$$M(n)=M^{\l}(n)\oplus W$$ 
for some suitable $\C^{\a_M}[G]$ submodule $W$ of $M(n).$ 
 
Define $\beta \in \End M(n)$ such that it is the identity on 
$M^{\l}(n)$ and zero on $W.$ Then $\beta\in (\End M(n))^{G}$. From the proof
of lemma \ref{corect} there  exists $v\in V^G $ such that 
$o_M(v)=\beta $ is the identity on $M^{\l}(n)$ and zero on $W.$ Thus, there is no $A_{n}(V^G)$-module 
homomorphism between $M_{\lambda}(n)$ and $M_{\gamma}(n).$ That is,
$M_{\lambda}(n)$ and $M_{\gamma}(n)$ are inequivalent $A_n(V^G)$-modules.

Since $M$ is irreducible $V$-module there exists $c\in \C$ such that
$M=\sum_{n\geq 0}M_{c+n}$ with $M_c\ne 0$ where $M_{c+n}$ is the
eigenspace for $L(0)$ with eigenvalue $c+n$ (see Definition \ref{d1}).
Thus we can take $M(n)=M_{c+n}.$ If $\dim (M_{\l})_{c+m}\ne \dim
(M_{\gamma})_{c+m}$ it is clear that $M_{\l}$ and $M_{\gamma}$ are nonisomorphic
$V^G$-module. Otherwise by Theorem \ref{l2.3}, 
$M_{\l}$ is not isomorphic to $M_{\gamma}$ as $V^G$-modules. \qed

\section{Dual pair II}
\setcounter{equation}{0}

Let $V$ be  a simple 
vertex operator algebra as in the last section 
and $G$ finite automorphism group of $V.$ Let $M$ be an irreducible
$V$-module. But we do not assume that $M$ is $G$-stable. 
We set $$G_M=\{g\in G| M\circ h\cong M\}$$
which is a subgroup of $G.$ Recall Theorem \ref{t2y}. One of
the main results in this section
is the following result of Schur-Weyl type:

\begin{th}\label{t3y} With the same notation as above we have:

1) $M^{\l}$ is nonzero for any $\l\in \Lambda_{G_M,\a_M}.$

2) Each $M_{\l}$ is an irreducible $V^G$-module.

3) $M_{\l}$ and $M_{\gamma}$ are equivalent $V^G$-module if and only
if $\l=\gamma.$

That is, $V^G$ and $\C^{\a_M}[G_M]$ form a dual pair on $M.$
\end{th}

This result generalizes Theorem 2.4 of [DLM1] (a case $M=V$) and sharpens
Theorem \ref{t2y}. In particular, the result shows that $M$ is a completely
reducible $V^G$-module. 

In order to prove Theorem \ref{t3y} we need a general setting
and we will prove a stronger result. 

\begin{de} A set $\S$ of irreducible $V$-modules is called {\em 
stable} if for any $M\in \S$ and $g\in G$ there exists $N\in \S$ such that
$M\circ g\cong N.$
\end{de}

Assume that ${\S}$ is a finite set of inequivalent irreducible $V$-module 
which is $G$-stable. Since
$M\circ(g_1g_2)$ and $(M\circ g_1)\circ g_2$ are isomorphic $V$-module
for $g_1,g_2\in G$ and $M\in \S$ we can 
define an right action of $G$ on $\S.$ Let $M\in \S$ and $g\in G$
we define $M\cdot g=N$ if $M\circ g\cong N.$ It is clear that
this action makes $\S$ a right $G$-set. 

\begin{rem} {\rm 1) For any irreducible $V$-module $M$ consider
the coset decomposition $G=\cup_{i=1}^kG_Mg_i.$ Then ${\S}=\{M\circ g_i|i=1,...,k\}$ is a such right $G$-set which will be used to prove Theorem  \ref{t3y}.

2) If $V$ is rational, a complete list
of inequivalent irreducible $V$-modules forms a such right $G$-set.}
\end{rem} 

Let $M\in {\S}$ and $x\in G$. Then 
there exists $N\in {\S}$ such that 
$N\cong M\circ x.$ That is,  
there is a linear map 
$\p_N(x):N\rightarrow M$ satisfying the condition:
$\p_N(x)Y_N(v,z)\p_N(x)^{-1}=Y_M(xv,z).$ 
By simplicity of $N$, there exists $\a_N(y,x)\in \C^*$ 
such that $\p_M(y)\p_N(x)=\a_N(y,x)\p_N(yx).$ 
Moreover, for $x,y,z\in G$ we have 
$$\a_N(z,yx)\a_N(y,x)=\a_M(z,y)\a_N(zy,x).$$ 

As in Section 3 we set $\C{\S}=\bigoplus_{M\in {\S}}\C e(M)$ 
and $e(M)e(N)=\d_{M,N}e(M)$. Let $U(\C{\cal S})
=\{\sum_{M\in{\S}}\l_{M}e(M)|\l_{M}\in \C^*\}$ and 
$\a(h,k)=\sum_{M\in {\S}}\a_{M}(h,k)e(M).$ 
It is easy to check that $\a(hk,l)\a(h,k)^l=\a(h,kl)\a(k,l).$ So 
$\a\in Z^2(G,U(\C\cal{S}))$ is a 2-cocycle.
Following Section \ref{s3}, we construct an associative algebra 
$A_{\a}(G,{\S})$ with multiplication defined by 
$$ a\o e(M)\cdot b\o e(N)=\a_{N}(a,b)ab\o e(M\cdot b)e(N)$$
for $a,b\in G$ and $M,N\in\S.$ 

We define an action of $A_{\a}(G,{\S})$
on $\oplus_{N\in \S}N$
in the following way: for $M,N\in {\S}$ and $w\in N$ we set
$$a\o e(M)\cdot w= \d_{M,N}\p_{M}(a)w$$
where $\p_{N}(a):N\rightarrow N\cdot a^{-1}.$  

For $N\in {\S}$ we let ${\O}_N=\{N\cdot g|g\in G\}$ 
be the orbit of $N$ under $G.$ 

\begin{lem} With the action defined above, $\oplus_{N\in \S}N$
becomes a module for $A_{\a}(G,{\S}).$ Moreover $\oplus_{M\in \O_N}M$ 
is a submodule for any $N.$
\end{lem}
\pf The proof is a straightforward computation: for $a,b\in G$ and
$M,N\in \S$ and $m\in N,$ 
\begin{eqnarray*}
& &\ \ (b\o e(L)\cdot a\o e(M))\cdot m\\
& &=(\d_{L\cdot a,M}\a_{M}(b,a)ba\o e(M))\cdot m\\
& &=\d_{L\cdot a, M}\d_{M,N}\a_{M}(b,a)\p_{M}(ba)m\\
& &=\d_{L\cdot a, M}\d_{M,N}\p_{M\cdot a^{-1}}(b)\p_{M}(a)m
\end{eqnarray*}
and
\begin{eqnarray*}
& &\ \ b\o e(L)\cdot (a\o e(M)\cdot m)\\
& &= b\o e(L)\cdot (\d_{M,N}\p_{M}(a)m)\\
& &=\d_{L,M\cdot a^{-1}}\d_{M,N}\p_{L}(b)\p_{M}(a)m\\
& &=\d_{L,M\cdot a^{-1}}\d_{M,N}\p_{M\cdot a^{-1}}(b)\p_{M}(a)m.
\end{eqnarray*}
Thus  $\bigoplus_{M\in \S}M$ is 
$A_{\a}(G,{\cal S})$-module. It is clear that  $\oplus_{M\in \O_N}M$ 
is a submodule.
\qed

Let $M\in \S.$  Recall from Section 3 that 
$D(M)=\<a\o e(M)|a\in G\>$ and 
$S(M)=\<a\o e(M)|a\in G_M\>$. 

\begin{prop}\label{p4.27} Let $N\in\S.$ Then $N$ is a $S(N)$-module
and there is an $A_{\a}(G,{\S})$-modules isomorphism between 
$D(N)\o_{S(N)}N$ and $\bigoplus_{M\in {\cal O}_N}M$ 
determined by $\Psi: a\o e(N)\o m \mapsto \p_{N}(a)m.$
\end{prop}
\pf To prove $N$ is a $S(N)$-module, it is enough to prove that 
$N$ is $S(N)$-stable. But this clear as $a\otimes e(N)w=\phi_N(a)w\in N$
for $a\in G_N$ and $w\in N.$ 

Next we prove that $\Psi$ is well defined. We must show that
$\Psi(d_{a,N}d_{b,N}\otimes w)= \Psi(d_{a,N}\otimes d_{b,N}w)$ 
for $w\in N,$ $a\in G, b\in G_N$
where $d_{a,N}=a\otimes e(N).$  
This
is clear as
\begin{eqnarray*}
& &\ \ \ \ \Psi(d_{a,N}d_{b,N}\otimes w)\\
& &=\alpha_N(a,b)\Psi(d_{ab,N}\otimes w)\\
& &=\alpha_N(a,b)\phi_N(ab)w\\
& &=\phi_N(a)\phi_N(b)w\\
& &=\phi_N(a)d_{b,N}w\\
& &=\Psi(d_{a,N}\otimes d_{b,N}w).
\end{eqnarray*}

The fact that $\Psi$ is a module homomorphism follows from
\begin{eqnarray*}
& &\Psi(d_{b,M}\cdot d_{a,N}\o w)\\
& &=\Psi(\d_{M\cdot a,N}\a_{N}(b,a)d_{ba,N}\o w)\\
& &=\d_{M\cdot a,N}\a_{N}(b,a)\p_{N}(ba)w)\\
& &=\d_{M\cdot a,N}\p_{N\cdot a^{-1}}(b)\p_{N}(a)w
\end{eqnarray*}
and
\begin{eqnarray*}
& &d_{b,M}\Psi(d_{a,N}\o w)=d_{b,M}\cdot \p_{N}(a)w\\
& &=\d_{M,N\cdot a^{-1}}\p_{N\cdot a^{-1}}(b)\p_{N}(a)w\\
\end{eqnarray*}
where $a,b\in G.$

In order to show that $\Psi$ is a bijection we construct an inverse
of $\Psi.$ 

Suppose that 
$ G=\bigcup_{i_N=1}^{k_N}G_Ng_{i_N}$ is a coset decomposition. Then 
$$\bigoplus_{M\in {\O}_N}M=\bigoplus_{i_N=1}^{k_N}
\p_N(g_{i_N}^{-1})N.$$
Define
\begin{eqnarray*}
\bigoplus_{M\in {\O}_N}M&\stackrel{\ch}{\rightarrow}& D(N)\bigotimes_{S(N)}N\\
\sum_{i_N=1}^{k_N}\p_{N}(g_{i_N}^{-1})m_{i_N}&\mapsto& \sum_{i_N=1}^{k_N}d_{g_{i_N}^{-1},N}\o m_{i_N}.
\end{eqnarray*}
It is straightforward to show that 
$\Psi\circ \ch=Id_{\bigoplus_{M\in {\O}_N}M}$ 
and $\ch\circ \Psi=Id_{D(N)\o_{S(N)}N}$. 
 So, ${D(N)\bigotimes_{S(N)}N}$ and $\bigoplus_{M\in {\O}_N}M$ are isomorphic as $A_{\a}(G,\cal {M})$-modules.\qed

Set $${\cal M}=\bigoplus_{M\in {\S}}M.$$
\begin{prop}\label{p50} 1) Every simple  
$A_{\a}(G,{\cal S})$-module occurs as a submodule of $\M$.

2) The actions of $A_{\a}(G,{\cal{M}})$ and $V^G$ on ${\cal M}$ commute.
\end{prop}

\pf Part 1) follows from Theorems \ref{t50},  \ref{t2y} and a natural identification
of $S(N)$ with $\C^{\a_N}[G_N]$ for
$N\in \S.$  Part 2) follows from
the relation 
$$a\otimes e(M)\bar Y(v,z)=\bar Y(av,z)a\otimes e(M)$$
on $\cal M$ for $a\in G,$ $M\in \S$ and $v\in V$ 
where we have used $\bar Y(v,z)$ denote the operator on $\cal M$ which acts
on $M$ by $Y_M(v,z).$ 
\qed

\begin{lem} Set $B=G\times U(\C{\S})$.  Then $B$ is a group under 
multiplication 
defined by $$(y,v)(x,u)=(yx,\a(y,x)v^xu)$$
\end{lem}

This result is standard as $\a$ is a 2-cocycle.

\begin{lem} $\M$ is a $B$-module via 
$(x,\sum_{M\in {\S}}\l_{M}e(M))=\sum_{M\in \S}\p_{M}(x)\l_M.$
\end{lem}

This result is immediate by noting that the subset
$B'=\{a\otimes u|a\in G,u\in U(\C\S)\}$ of $A_{\alpha}(G,\S)$ 
is a multiplicative group 
isomorphic to $B$ and the action of $B$ on $\M$ is
exactly the action of $B'$ on ${\cal M}.$  

\begin{lem}
For any $(x,u)\in B$, $v\in V,$ we have 
$$(x,u)\bar{Y}(v,z)(x,u)^{-1}=\bar Y(xv,z)$$
on $\cal M.$
\end{lem}
\pf Let $b=(x,\sum_{M\in{\S}}\l_{M}e(M)))\in B$. 
Then
\begin{eqnarray*}
& &b\bar Y(v,z)=(x,\sum_{M\in{\S}}\l_{M}e(M)))\bar Y(v,z)\\
& &\ \ \ \ =\bar Y(xv,z)(x,\sum_{M\in{\S}}\l_{M}e(M)))
\end{eqnarray*}
(see the proof of Proposition \ref{p50} 2)) as required.
\qed

\begin{lem}\label{l50}
$\bigoplus_{M\in \S}\End(M)$ is a $B$-module via 
$$(x,\sum_{M\in{\S}}\l_{M}e(M))\cdot \sum_{M\in \S}f_M
=\sum_{M\in \S}\p_{M}(x)f_{M}\p_{M}(x)^{-1}.$$ 
Furthermore, for each orbit $\O$ and $n\geq 0,$ 
$\bigoplus_{M\in {\O}}\End(M(n))$ 
is a $B$-submodule. In particular, $\bigoplus_{M\in {\O}}\End(M(n))$ 
is a $G$-module.
\end{lem}

\pf 
Let $(x,\sum_{M\in{\S}}\l_{M}e(M))$, 
$(y,\sum_{M\in {\S}}\b_{M}e(M))\in B.$ Since 
\begin{eqnarray*}
& &[(y,\sum_{M\in {\S}}\b_{M}e(M))(x,\sum_{M\in{\S}}\l_{M}e(M))]\sum_{M\in {\S}}f_M    \\
& &=(yx,\sum_{M\in {\S}}\b_{M\cdot x^{-1}}\l_{M}\a_{M}(y,x)e(M))\sum_{M\in {\S}}f_M\\
& &=\sum_{M\in {\S}}\p_{M}(yx)f_{M}\p_{M}(yx)^{-1}\\ 
& &=\sum_{M\in {\S}}\p_{M\cdot x^{-1}}(y)\p_{M}(x)f_{M}\p_{M}(x)^{-1}\p_{M\cdot x^{-1}}(y)^{-1}
\end{eqnarray*}
and
\begin{eqnarray*}
& &(y,\sum_{M\in {\S}}\b_{M}e(M)) [(x,\sum_{M\in{\S}}\l_{M}e(M))\sum_{M\in {\S}}f_{M}]\\
& &(y,\sum_{M\in {\S}}\b_{M}e(M))\sum_{M\in {\S}}\p_{M}(x)f_{M}\p_{M}(x)^{-1}\\
& &=\sum_{M\in {\S}}\p_{M\cdot x^{-1}}(y)\p_{M}(x)f_{M}\p_{M}(x)^{-1}\p_{M\cdot x^{-1}}(y)^{-1},
\end{eqnarray*}
$\bigoplus_{M\in {\S}}\End(M)$ is a $B$-module. The other assertions
in the lemma are clear. 
\qed

\begin{rem}{\rm It is worth to point out that $\phi_M(x)f_M\phi_M(x)^{-1}$
does not lie in $\End M$ for $f_M\in \End M$ and $x\in G$ unless $x\in G_M.$
In general,  $\phi_M(x)f_M\phi_M(x)^{-1}$ is an element of $\End (M\cdot x^{-1}).$}
\end{rem}

For any nonnegative $n\in\Z,$ let $\sigma_n$ be a map from $A_n(V)$ to $\bigoplus_{M\in \S}\sum_{0\leq m\leq n}\End M(m)$ 
defined by $\sigma_n(v+O_n(V))=\sum_{M\in \S}o_{M}(v).$
\begin{lem}\label{l51}
The map $\sigma_n$ is a $G$-module epimorphism. In particular,
$\sigma(A_n(V)^G)=(\sum_{M\in{\S},0\leq m\leq n}\End M(m))^G$ 
\end{lem}
\pf  Let $g\in G$ and $v\in V.$ 
Then 
\begin{eqnarray*}
& &\sigma_n(g(v+O_n(v)))=\sigma_n(gv+O_n(V))\\
& &\ \ \ \ \ =\sum_{M\in \S}o_{M}(gv)\\
& &\ \ \ \ \ =\sum_{M\in \S}\p_M(g)o_M(v)\p_M(g)^{-1}.
\end{eqnarray*}
That is, $\sigma_n$ is a $G$-homomorphism. 

In order to see that $\sigma_n$ is onto we note that all $M\in \S$ are
inequivalent. We assume that $M(0)\ne 0$ for all $M\in \S.$ 
By Theorem \ref{l2.3}, $\{M(m)|M\in {\S}, 0\leq m\leq n\}$ is a set
of finite dimensional inequivalent $A_n(V)$-modules. Let $K_{M(m)}$ be
the kernel of $A_n(V)$ on $M(m).$ Then $A_n(V)/K_{M(m)}$ is
isomorphic to $\End M(m)$ and $A_n(V)/K_n$ is isomorphic to
the direct sum $\bigoplus_{M\in \S}\bigoplus_{0\leq m\leq n}\End M(m)$ where $K_n=\cap_{M,m}K_{M(m)}.$ Thus $\sigma_n$ is onto.
\qed

Now we are in the position to state and to prove the main result in this paper.
Let ${\S}=\cup_{j\in J}\O_j$ be an orbit decomposition and
fix $M^j\in \O_j.$ For convenience, we set $G_j=G_{M^j}$ and
$\Lambda_j=\Lambda_{M^j,\alpha_{M^j}}.$ 
Then by Theorem \ref{t2y} we have
a decomposition
$$M^j=\sum_{\l\in \L_j}W_{\l}\otimes M^j_{\l}$$
where $M^j_{\l}$ is an irreducible $V^{G_j}$-module.
Thanks to Propositions \ref{p4.27} and \ref{p50}  we have
$${\M}=\bigoplus_{j\in J}\sum_{\l\in {\L}_j}(\Ind_{S(M^j)}^{D(M^j)}
W_{\l})\otimes M^j_{\l}$$ 
as a $A_{\alpha}(G,{\S})\otimes V^G$-module. For $j\in J$ and
$\l\in {\L}_j$ we set $W_{j,\l}=\Ind_{S(M^j)}^{D(M^j)}
W_{\l}.$ Then by Theorems \ref{t50}, $W_{j,\l}$ forms a complete
list of simple $A_{\alpha}(G,\S)$-modules. 

\begin{th}\label{main} As a $A_{\alpha}(G,{\S})\otimes V^G$-module,
$${\M}=\bigoplus_{j\in J,\l\in \L_j}W_{j,\l}\otimes M^j_{\l}.$$
Moreover,

1) Each $M^j_{\l}$ is a nonzero irreducible $V^G$-module.

2) $M^{j_1}_{\l_1}$ and $M^{j_2}_{\l_2}$ are isomorphic $V^G$-module
if and only if $j_1=j_2$ and $\l_1=\l_2.$ 

That is, $(A_{\alpha}(G,{\S}), V^G)$ forms a dual pair on $\M.$
\end{th}

\pf 1) We have already seen from Theorem \ref{t2y} that each $M^j_{\l}$
is nonzero. Again by Theorem \ref{l2.3} 3) we only need to
show that each $M^j_{\l}(n)=M^j(n)\cap M^j_{\l}$ is a simple 
$A_n(V^G)$-module. 

We now fix $j.$ Suppose $G=\bigcup_{i=1}^{k}G_{j}g_{i}$ be the coset 
decomposition
such that $g_1=1.$ 
For each $v\in V^{G_j}$ we set $f^v\in \sum_{M\in \S}\End M(n)$
such that $f^v$ acts on $M^j\cdot g_i^{-1}$ as 
$$\phi_{M^j}(g_i)o_{M^j}(v)
\phi_{M^j}(g_i)^{-1}=o_{M^j\cdot g_i^{-1}}(g_iv)$$ 
for $i=1,...,k,$ acts on any other $M(n)$
as zero. Then it is clear from Lemma \ref{l50} that $f^v\in 
(\sum_{M\in \S}\End M(n))^G.$ By Lemma \ref{l51} there exists
$u\in V^G$ such that $f^u=f^v.$ This shows that for any
$v\in V^{G_{j}}$ there exists $u\in V^G$ such that $o_{M^j}(v)=o_{M^j}(u)$
on $M^j(n).$ Since $M^j_{\l}(n)$ is a simple module for
$A_n(V^{G_{j}})$ by Theorems \ref{l2.3} and \ref{t2y}
we see immediately that $M^j_{\l}(n)$ is a simple $A_n(V^G)$-module.

2) We take $(j_1,\l_1)\ne (j_2,\l_2).$ We must to prove that
 $M^{j_1}_{\l_1}$ and $M^{j_2}_{\l_2}$ are inequivalent.
Let $n\geq 0$ such that both $M^{j_1}_{\l_1}(n)$ and 
$M^{j_2}_{\l_2}(n)$ are nonzero. 
Then $$\sum_{0\leq m\leq n}{\M}(m)
=\sum_{0\leq m\leq n}W_{j_1,\l_1}\otimes M^{j_1}_{\l_1}(m)\oplus W$$ for a
suitable $A_{\alpha}(G,\S)$-module $W.$ 
Define $f$ in $\sum_{M\in {\S},0\leq m\leq n}\End M(m)$ such that $f=1$ on
$\sum_{0\leq m\leq n}W_{j_1,\l_1}\otimes M^{j_1}_{\l_1}(m)$ and $f=0$ on $W.$ 
Again it is obvious that 
$$f\in (\sum_{M\in {\S},0\leq m\leq n}\End M(m))^G.$$
Using Lemma \ref{l51} we find $v\in V^G$ such that
$\sigma_n(v)=f.$ That is $o_{M^{j_1}_{\l_1}}(v)=1$
on $\sum_{0\leq m\leq n}M^{j_1}_{\l_1}(m)$ and $o_{M^{j_2}_{\l_2}}(v)=0$
 on $\sum_{0\leq m\leq n}M^{j_2}_{\l_2}(m).$ 
That is,  $M^{j_1}_{\l_1}(s)$ and $M^{j_2}_{\l_2}(t)$
are nonisomorphic $A_n(V^G)$-modules for $0\leq s,t\leq n$ if either
 $M^{j_1}_{\l_1}(s)$ or $M^{j_2}_{\l_2}(t)$ is nonzero. 
In particular, $M^{j_1}_{\l_1}(s_0)$ or $M^{j_2}_{\l_2}(t_0)$
are nonisomorphic $A_n(V^G)$-modules where
$s_0,t_0\geq 0$ such that   $M^{j_1}_{\l_1}(s)=0$
and $M^{j_2}_{\l_2}(t)=0$ for all $s\leq s_0$ and $t\leq t_0.$ 
Our choices of $s_0$ and $t_0$ then assert that 
$M^{j_1}_{\l_1}(s_0)$ or $M^{j_2}_{\l_2}(t_0)$ are, in fact, inequivalent
$A_0(V^G)$-modules. 
Thus by Theorem \ref{l2.3} 
$M^{j_1}_{\l_1}$ and $M^{j_2}_{\l_2}$ are inequivalent
$V^G$-modules. \qed

Now Theorem \ref{t3y} follows from Theorem \ref{main} immediately by
taking the right $G$-set $\S$ to be the $G$-orbit
$\{M\circ g_i|i=1,...,l\}$ where $G=\cup_{i=1}^l G_Mg_i$
is the coset decomposition of $G$ with respect to the stabilizer
$G_M=\{g\in G|M\circ g\cong M\}.$

We end this paper with the following general discussion:  If
$V$ has only finitely many inequivalent irreducible modules
(this is the case when $V$ is rational; see the discussion after Definition 
\ref{ration}), then 
a complete list of irreducible 
$V$-modules is a right $G$-set. Theorem
\ref{main} then tells us not only ${\M}=\sum_{M\in \S}M$
is completely reducible but also gives an equivalence
between $A_{\alpha}(G,\S)$-module category 
and a subcategory of admissible $V^G$-modules 
generated by the irreducible submodules occurring in
$\M$ by sending  $W_{j,\l}$ to $M^j_{\l}.$ We expect to
use this result to determine the module category for
$V^G$ when $V$ is rational in the future.

\end{document}